\documentclass{eptcs}
\providecommand{\event}{ACT2022}

\title{Fibrational Linguistics (\fiblang): Language Acquisition}
\makeatletter%
\@ifclassloaded{amsart}{
  \usepackage{fouche/fouche}
\author{Fabrizio \textsc{Genovese}}
\author{Fosco \textsc{Loregian}}
\author{Caterina \textsc{Puca}}
\address{ Fabrizio Genovese\newline
	University of Pisa / Statebox,\newline %
	\url{fabrizio.romano.genovese@gmail.com}%
}
\address{ Fosco Loregian\newline
	Tallinn University of Technology,\newline %
	\url{fosco.loregian@taltech.ee}%
}
\address{ Caterina Puca\newline
	Sapienza University of Rome,\newline %
	\url{caterpuca@gmail.com}%
}
\usepackage{hyperref}
\hypersetup{ hyperfootnotes=true
	, bookmarksnumbered=true
	, bookmarksopen=true
	, breaklinks=true
	, urlcolor=red!40!black
	, colorlinks%
	, urlcolor=RoyalBlue
	, linkcolor=Red
	, citecolor=Green
}
\usepackage{tipa}

\thanks{The first author was supported by the \href{https://gitcoin.co/grants/1086/independent-ethvestigator-program}{Independent Ethvestigator Program}. The second author was supported by the ESF funded Estonian IT Academy research measure (project 2014--2020.4.05.19--0001). A special thanks goes to Bettino, the fat cat of figures 1, 2, and 3, and Peppa the black evil feline, for having provided valuable suggestions and
comments after an attentive proofreading of an early version of this manuscript.}
}%
{
  \usepackage{hyperref}
\hypersetup{ hyperfootnotes=true
	, bookmarksnumbered=true
	, bookmarksopen=true
	, breaklinks=true
	, urlcolor=red!40!black
	, colorlinks%
	, urlcolor=RoyalBlue
	, linkcolor=RoyalBlue
	, citecolor=Green
}
\usepackage{fouche/act-fouche}
\author{ Fabrizio \textsc{Genovese}
	\email{\orcidNum{0000-0001-7792-1375}}
	\institute{20[]}%
	  \thanks{Supported by the \href{https://gitcoin.co/grants/1086/independent-ethvestigator-program}{Independent Ethvestigator Program}.}
	\email{fabrizio.romano.genovese@gmail.com}
	\and
	Fosco \textsc{Loregian}
	\email{\orcidNum{0000-0003-3052-465X}}
	\institute{Tallinn University of Technology}%
	  \thanks{Supported by the ESF funded Estonian IT Academy research measure (project 2014--2020.4.05.19--0001).}% A special thanks goes to Bettino, the fat cat of figures 1, 2, and 3, and Peppa the black evil feline, for having provided valuable suggestions and	comments after an attentive proofreading of an early version of this manuscript.}
	\email{fosco.loregian@gmail.com}
	\and
	Caterina \textsc{Puca}
	\institute{Quantinuum, 17 Beaumont Street, Oxford, OX1 2NA, United Kingdom}
	\email{caterpuca@gmail.com}
}
\def\authorrunning{Genovese%
	, Loregian%
	, Puca%
}
\def\titlerunning{Fibrational linguistics}

\usepackage{ newtxtext
         	 , newtxmath
           }

}%
\makeatother%
\usetikzlibrary{decorations.pathreplacing,calligraphy}

\def\Gro#1{{\textstyle\int} #1}
\def\GroI#1{\nabla #1}
\def\fiblang{FibLang\@\xspace}
\def\catIcon{C}
\def\catIconSmall{C}

\newcommand{\framew}[1]{{{#1}^\sharp}}

\usepackage{doi}

\usepackage{doi}
\usepackage[ backend=bibtex
	         , natbib=true
	         , style=numeric
	         , sorting=nyt
	         , maxbibnames=10
	         , giveninits=true
	         , backref=true
	         , url=false
	         , doi=true
	         , isbn=false
	         , arxiv=abs
	         , date=year
	         , abbreviate = true
           ]{biblatex}

\AtBeginBibliography{\small}

\def\lCell{\xymatrix{
\ar@{>->}[d]_{e_f} \ar@{}[dr]|{\alpha_L} \ar[r]|*@{|}^p& \ar@{>->}[d]^{e_g}\\
\ar[r]|*@{|}_q&
}}
\def\rCell{\xymatrix{
\ar@{->>}[d]_{m_f} \ar@{}[dr]|{\alpha_R} \ar[r]|*@{|}^p& \ar@{->>}[d]^{m_g}\\
\ar[r]|*@{|}_q&
}}
\def\Cell{\xymatrix{
\ar[d]_f \ar@{}[dr]|\alpha \ar[r]|*@{|}^p& \ar[d]^g\\
\ar[r]|*@{|}_q&
}}

\def\cap#1{\todo[inline, caption=0, color=blue!10]{Cat said: #1}}

\usepackage{xspace}
\usepackage{graphbox}

\usepackage[english]{babel}

\newcommand{\orcidLogo}{\includegraphics[align=c]{orcid_logo.pdf}\kern.3em}
\newcommand{\orcidNum}[1]{\href{https://orcid.org/#1}{#1}}

\NewDocumentCommand{\fib}{O{r}}{\ar@{->>}[#1]|*@{*}}
\NewDocumentCommand{\opfib}{O{r}}{\ar@{->>}[#1]|*@{o}}

\def\Gro#1{{\textstyle\int} #1}
\def\GroI#1{\nabla #1}
\def\fiblang{FibLang\@\xspace}
\def\catIcon{C}
\def\catIconSmall{C}

\usepackage{cleveref}

\usepackage{dialogue,csquotes}
\usepackage{anyfontsize, wrapfig}
\usepackage{tempora}
\usepackage[super]{nth}

\usetikzlibrary{decorations.pathmorphing}

\begin{document}

\maketitle
\begin{abstract}
  In this work we show how \fiblang, a category-theoretic framework concerned with the interplay between language and meaning, can be used to
  describe \emph{vocabulary acquisition}, that is the process with which a speaker $p$ acquires new vocabulary (through experience or interaction).

  We model two different kinds of vocabulary acquisition, which we call `by example' and `by paraphrasis'. The former captures the idea of acquiring the meaning of a word by being shown a witness representing that word, as in `understanding what a cat is, by looking at a cat'. The latter captures the idea of acquiring meaning by listening to some other speaker rephrasing the word with others already known to the learner.

  We provide a category-theoretic model for vocabulary acquisition by paraphrasis based on the construction of free promonads. We draw parallels between our work and Wittgenstein's dynamical approach to language, commonly known as 'language games'.
\end{abstract}

\section{Introduction}
Language has always been characterised as a distinctive, exclusive feature of human beings, yet children are not \emph{born} fluent in any language at all. Along the history of human thought, this apparent discrepancy has promptly led to philosophical speculation on the innateness of language \cite{plato1949meno}, which was subsequently replaced with a more cautious theory of innateness of syntactic structures \cite{chomsky2006language}. On the other hand, empiricists such as J. Locke \cite{locke1847essay} argued that the human mind had to be thought of as a \emph{tabula rasa}.

Either way, these heterogeneous philosophical stances share the necessity of formalising a common process: \emph{language acquisition}, i.e. the process in which proficiency in a language increases with time or solicitation.\footnote{Here and in the rest of paper we will be referring to first language acquisition only: although similarities have been pointed out concerning second language acquisition, such as \emph{the silent period} \cite{ellis1989understanding}, several substantial differences separate the two processes \cite{ipek2009comparing}. For instance, there is evidence that the learner's first language slows the development of acquisitional sequences predicted by the Natural Order Hypothesis \cite{krashen2002theory, mclaughlin1987theories}. Additionally, according to the Critical Period Hypothesis, after puberty, lateralisation is accomplished, and reduced plasticity of the brain can compromise the fluent acquisition of a second language \cite{brown2000principles}.} Although years of debate in linguistics have not come to a definite resolution yet, the transformationalist orientation of Chomsky has received severe criticism in the \nth{20} century with the growing development of linguistic philosophy and with the renewed interest in L. Wittgenstein's ideas.

In particular, Chomsky's model of language acquisition has been accused of being overly reductionist and mechanical \cite{phillips1971ludwig}, as opposed to Wittgenstein's dynamic theory that embraces \emph{context} as a fundamental aspect of meaning analysis \cite{Wittgenstein2009, wittgenstein1958blue}.

In this latter perspective, syntax and semantics must interact and reciprocally influence each other during communication, a process also referred to as a \emph{language game} \cite{chakraborty2018language}. More precisely, language games can serve as a tool to untie the problems of context-dependency and ambiguity regarding words with multiple semantic interpretations.

The recently developed \emph{FibLang} framework~\cite{fiblang0} takes a stab at tackling the enticing and deep problem of language, offering a category\hyp{}theoretic framework concerned with describing the interplay between meaning and structure in natural language. As a theory, \fiblang relies on fibered categories \cite{CLTT,loregian2019categorical,Vistoli2005, street1980fibrations}; the main idea underlying \fiblang is characterizing linguistic meaning as \emph{fibered over grammar}.

Here we argue that Wittgenstein's perspective hints at a fibrational formalisation of language acquisition. We will substantiate our hypothesis by formalising vocabulary acquisition in \fiblang and show how the context-dependency and semantic ambiguity aspects underlying language games are organically embraced in our categorical description.

Since our fibrational approach to language acquisition naturally encodes agency, it is out-of-the-box compatible with applications to learning tasks in natural language processing (NLP). In this sense, it also enriches the static perspective of DisCoCat \cite{coecke2010mathematical}, which by the way, has also been recently revised using tools from categorical game theory applied to language games in \cite{8bd9537e5b18e941c46c313c3feb4edbc542e345}.

In point of fact, in more recent years, a different framework, named DisCoCirc, has been adopted to allow for a dynamic flow between syntax and semantics \cite{coecke2021compositionality}. Implementation aspects regarding quantum computers suggested this switch and, as a side-effect, it provides stronger foundations to the philosophical stance of \fiblang.

\paragraph*{Structure of the paper.}
In~\autoref{sec: rappels}, we will recall the basic definitions of FibLang. We will then provide our description of vocabulary acquisition by example in \autoref{sec: vocabulary acquisition by example}. Subsequently, we will define the tool of explanations in~\autoref{sec: explanations}, and will use them to define vocabulary acquisition by paraphrasis in \autoref{sec: vocabulary acquisition by paraphrasis}. Finally, at the end of \autoref{sec: vocabulary acquisition by paraphrasis}, we provide a construction (cf. \autoref{shit_construction}) with which to show how our formalisation of vocabulary acquisition can be used to enrich grammar by appropriately acknowledging semantic interrelations.
\section{Rappels of \fiblang}\label{sec: rappels}
\fiblang was introduced in a prior installment of this series, \cite{fiblang0}. Its main idea can be summarized as follows: whereas it is reasonable to believe that language has at least some degree of compositionality, especially when describing grammar, it becomes much more difficult to substantiate this position when it comes to describe meaning. Indeed, compositional --and especially cognitive-- models of meaning such as G\"ardenfors' \cite{gardenfors2004conceptual} are prone to criticism on multiple fronts, not last the fact that a truly universal model of meaning is very difficult to define because of cultural and cognitive differences between speakers.

To circumvent these problems \fiblang focuses on describing the interplay between meaning and structure in abstract terms in a way that is agnostic to the particular model one chooses to represent either. This vision is reified in the main idea being that a `language' is a category $\clL$ of some sort, and a speaker $p$ of a language $\clL$ is a \emph{fibration} $\var[\framew{p}]{\clE^p}{\clL}$ over $\clL$.\footnote{In \cite{fiblang0} fibrations are denoted using the superscript $\framew{-}$ to distinguish them visually from bare functors.} The defining property of a fibration is that its domain $\clE^p$ is obtained gluing together all the various fibres $\clE_L$ over objects of $\clL$ in a coherent manner.

Borrowing from the work of Lambek \cite{lambek1958mathematics}, $\clL$ represents a `language category' while borrowing from fibration theory \cite{benabou1985fibered,CLLT} $\clE^p$, called the \emph{total category} of the fibration, represents a `semantics category' for the speaker $p$. This semantics category could be thought of as any sort of cognitive or distributional model of meaning for the language in question: this is to say that we are not particularly attached to any specific model but rather aim at the highest possible generality.

To arrive at this idea, the starting point in \cite{fiblang0} constitutes the most general and yet reasonable assumption one could make, namely \emph{there is some structure-preserving map from what a speaker means to what a speaker says}. Formally, this directly translates to modelling speakers as simple functors between categories.
From this, multiple reasons are stated that lead to believe that the language category $\clL$ should be treated as something that can be explicitly modelled and studied, while $\clE^p$ should be treated as a black box. Then, it is shown how every functor can be factorised into a fibration via~\autoref{thm:_factorization}, obtaining a more workable definition of a speaker from the very abstract one we started from.
\begin{definition}[Fibration]\label{def:_fibration}
  A functor $\framew{p} : \clE \to \clC$ is a \emph{discrete fibration} (for us, just a \emph{fibration}) if, for every object $E$ in $\clE$ and every morphism $f: C \to \framew{p}E$ with $C$ in $\clC$, there exists a unique morphism $h: E' \to E$ such that $\framew{p}h = f$.
\end{definition}
\begin{notation}\label{notat_for_fib}
  The domain of a fibration $\framew{p} : \clE \to \clC$ is usually called the \emph{total category} of the fibration, and its codomain is the \emph{base category}. Given any functor $p$ we can define the \emph{fibre} of $p$ over an object $C\in\clC$, i.e. the subcategory $\clE_C = \{f : E\to E'\mid p f = 1_C\}\subseteq \clE$.
\end{notation}
Intuitively, a fibration is a functor $\framew{p}: \clE \to \clC$ that realises the category $\clE$ as a `covering' of $\clC$, in such a way that morphisms in $\clC$ can be lifted to $\clE$, to induce functions between the fibres in $\clE$, called \emph{reindexing functions}. We have represented an elementary example of fibration in the following figure, with the action of the reindexing functions between the elements in the fibres made explicit. The gray rectangles are the \emph{fibres} of the fibration.
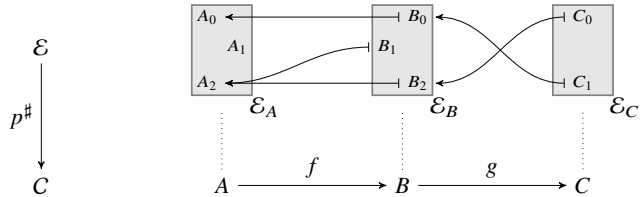
\begin{equation*}
  \scalebox{0.8}{
    \begin{tikzpicture}[node distance=1.3cm,>=stealth',bend angle=45,auto]
      \node (C) at (-3,0) {$\clC$};
      \node (E) at (-3,2.3) {$\clE$};

      \node (1a) at (0,0) {$A$};
      \node (3a) at (3,0)  {$B$};
      \node (5a) at (6,0) {$C$};

      \draw[->] (1a) to node[above] {$f$} (3a);
      \draw[->] (3a) to node[above] {$g$} (5a);

      \node (EA) at (0.7,1.3) {$\clE_A$};
      \draw[thick, draw=black!50, fill=black!10] (-0.5,1.5) rectangle (0.5,3);
      \node[font=\scriptsize] (a0) at (-.25,2.8) {$A_0$};
      \node[font=\scriptsize] (a1) at (0.25,2.3) {$A_1$};
      \node[font=\scriptsize] (a2) at (-.25,1.7) {$A_2$};

      \node (EB) at (3.7,1.3) {$\clE_B$};
      \draw[thick, draw=black!50, fill=black!10] (2.5,1.5) rectangle (3.5,3);
      \node[font=\scriptsize] (b0) at (3.25,2.8) {$B_0$};
      \node[font=\scriptsize] (b1) at (2.75,2.3) {$B_1$};
      \node[font=\scriptsize] (b2) at (3.25,1.7) {$B_2$};

      \node (EC) at (6.7,1.3) {$\clE_C$};
      \draw[thick, draw=black!50, fill=black!10] (5.5,1.5) rectangle (6.5,3);
      \node[font=\scriptsize] (c0) at (6,2.8) {$C_0$};
      \node[font=\scriptsize] (c1) at (6,1.7) {$C_1$};

      \draw[->] (E) to node[left] {$\framew{p}$} (C);

      \draw[<-|, out=0, in=180] (a0) to (b0);
      \draw[<-|, out=0, in=180] (a2) to (b1);
      \draw[<-|, out=0, in=180] (a2) to (b2);

      \draw[<-|, out=0, in=180] (b0) to (c1);
      \draw[<-|, out=0, in=180] (b2) to (c0);

      \draw[dotted, -]  (1a.north) -- (0,1.25);
      \draw[dotted, -]  (3a.north) -- (3,1.25);
      \draw[dotted, -]  (5a.north) -- (6,1.25);
    \end{tikzpicture}
  }
\end{equation*}
A fundamental result in the theory of fibrations that we will often use is that fibrations over a category $\clL$ are equivalent to functors out of $\clL$, with codomain the category of sets and functions. More details on the construction, and a full explanation of its usefulness for \fiblang, can be found in \cite[A.6]{fiblang0}; a classical reference for the theory of fibrations is \cite{CLTT}.
\begin{theorem}\label{thm: fibrations functors equivalence}
  There is a category $\cate{DFib}/\clL$ of fibrations over a given $\clL$, where an object is a fibration $\var[\framew{p}]{\clE^p}{\clL}$ and a map $h : \var[\framew{p}]{\clE^p}{\clL} \to \var[\framew{q}]{\clE^q}{\clL}$ is a functor $h : \clE^p \to \clE^q$ such that $\framew{q}\cdot h=\framew{p}$.

  There is an equivalence of categories:
  \[
    \GroI{-} : \cate{DFib}/\clL \cong [\clL^\op, \Set] : \Gro{-}
  \]
  where at the right-hand side we have the category of all functors $\clL\to\Set$ and natural transformations thereof. The functor $\Gro{-}$ is often called \emph{the category of elements construction}, or in its most general form \emph{the Grothendieck constuction}.
\end{theorem}
\autoref{thm: fibrations functors equivalence} is also instrumental in pointing out how \fiblang -- which postulates an approach going from meaning to language -- can be made compatible with traditional models of meaning such as DisCoCat -- which postulate an approach going from language to meaning\footnote{
  In the particular case of DisCoCat it would be more proper to say that the chosen approach is from \emph{grammar} to \emph{semantics}. The apparent dissonance is resolved by taking into account the agnostic approach of \fiblang. The language category $\clL$ can be purely grammatical, employing for instance a Lambek's pregroup~\cite{093fd03b7dc604c64124630f20a8f01231253397} as in DisCoCat, or more expressive. Similarly, the meaning category $\clE^p$ for a speaker $p$ could be purely semantical -- for instance, a distributional model or a conceptual space~\cite{gardenfors2004conceptual} -- or more expressive.
}. Many conceptual reasons are given in \cite{fiblang0} to prefer the meaning-to-language approach, but~\autoref{thm: fibrations functors equivalence} shows how the two are faces of the same medal.

As we remarked above, FibLang relies on some machinery to turn a model for speakers consisting of simple functors into a model consisting of fibrations. The main theorem allowing us to do so is the following:
\begin{theorem}[\protect{\cite[Theorem 3]{bams1183534973}}]\label{thm:_factorization}
  Any functor $p : \clD^p \to \clL$ can be written as a composition of functors $\clD^p \xrightarrow{s} \clE^p \xrightarrow{\framew{p}} \clL$, such that $\framew{p}$ is a fibration.
\end{theorem}
We will make heavy use of~\autoref{thm:_factorization} in the following sections to model language acquisition.
\section{Vocabulary acquisition by direct example}\label{sec: vocabulary acquisition by example}
\emph{Vocabulary acquisition} denotes the act of acquiring meanings for a word previously unknown \cite{mckeown2014nature}. In this work, we aim to describe two main modes of vocabulary acquisition. In this section, we focus on \emph{vocabulary acquisition by direct example}: this is the easiest method of language acquisition that we can describe and a commonly used method in children's education and monolingual fieldwork~\cite{everett2001monolingual,newman2001linguistic}: it simply works by pointing at something and saying the word one is referring to. This process of language acquisition, as Wittgenstein explains in \cite{wittgenstein2010philosophical}, can serve as a primitive example of a language game.
\begin{example}Consider the following dialogue:
    \begin{dialogue}\label{ex:dialogue}
        \speak{Alice} Look, a cat!
        \speak{Bob} A what?

        \direct{\refer{Alice} points to a cat}
        \speak{Alice} That, a cat!
        \speak{Bob} Oh!
    \end{dialogue}
    What happened in that `Oh!' can be mathematically modelled as a colimit in the fibrations that represent Alice and Bob.
\end{example}
\begin{definition}[Vocabulary acquisition by language example]\label{def:vocabulary acquisition by example}
    Consider two speakers $\var[\framew{p}]{\clE^p}{\clL}$ and $\var[\framew{q}]{\clE^q}{\clL}$, which we will call \emph{teacher} and \emph{learner}, respectively.

    Suppose that, for some $L \in \clL$ -- called \emph{the linguistic element to learn}\footnote{We will use the wording \emph{linguistic elements} referring to words, entire sentences or something else, depending on the model we chose for $\clL$, without committing to a particular choice. More formally, a linguistic element is the (possibly nonfull) subcategory of $\clL$ spanned by a certain choice of objects.} -- we have that $\clE^p_L \neq \emptyset$ and $\clE^q_L = \emptyset$.
    Fix a subset $S \subseteq \clE^p_L$, called an \emph{example} for $L$. Then we can define a new category $\clF^q$ as follows:
    \begin{align*}
        \text{obj}( {\clF^q} ) & := \text{obj}(\clE^q) \sqcup S \\
        \hom({\clF^q})            & := \hom(\clE^{q})
    \end{align*}
    and a functor $T : \clF^q \to \clL$ agreeing with $\var[\framew{q}]{\clE^q}{\clL}$ on every fibre $L'\neq L$, and sending every object of $S$ to $L$. Relying on~\autoref{thm:_factorization}, the new fibration modelling the speaker $q$ after learning $L$ is the factorization $\framew{\tilde{q}}$ such that:
    \begin{equation*}
        \xymatrix{
            T = \big(\clF^q \ar[r]^{s} & {\clE}^{\tilde{q}} \ar@{->>}[r]^{\framew{\tilde{q}}} & \clL\big).
        }
    \end{equation*}
\end{definition}
Let us unpack this definition. We consider two speakers $p,q$ of the same language $\clL$.
Speaker $q$ does not know the meaning of a given linguistic element which corresponds to an object $L \in \clL$ --which is why we call $q$ \emph{learner}. The fibre $\clE_L$ in $\clE$ is the empty set, as $L$ has no meaning for $q$.
On the other hand, speaker $p$ has some model of meaning for $L$ --which is why we call $p$ \emph{teacher}, and we assume the fibre $\clE^p_L$ to be not empty.

If $p$ points out an instance of $L$ to $q$, as in~\autoref{ex:dialogue}, it is reasonable to assume that the instance in question is itself part of the fibre $\clE^p_L$, as $p$ recognises the example as an instance of the concept $L$. So we postulate that the example identifies a subset $S \subseteq \clE^p_L$. While following along with the example, $q$ incorporates the set $S$ as a fibre over $L$ by extending its meaning.

On the intuitive side, there is no reason to claim that forcefully adding a scattered set of notions to the ones previously mastered by a speaker is enough to let $q$ `understand' that said set of notions defines a new language term. Instead, to attain this level of understanding, $q$ has to build meaningful relations between the new term $L$ and all the others they master, in concordance with all the pre-existing relations between general terms.

On the mathematical side, there is no reason why the functor obtained by forcefully adding a nonempty fibre to a previously empty one shall remain a fibration; to be such, we must let the new fibre interact well with the environment, in concordance with the hom-sets $\clL(X,L)$ and $\clL(L,Y)$. This is why we consider the comprehensive factorisation of $T$ and take into consideration only $\framew{\tilde{q}}$ instead of the whole $T$.
\def\w#1{\textsf{#1}}
\begin{example}[A language game]\label{exa:language game}
    Consider the following language game borrowed from Wittgenstein's \emph{Philosophical Investigations} \cite{wittgenstein2010philosophical}: a builder $A$ asks his assistant $B$ to pass him the stones with which they are building, in the order $A$ calls them out. In this situation, let us imagine that the language they use consists of only four words: \w{block}, \w{pillar}, \w{slab}, \w{beam}.

    This language game can be interpreted as a particular case of vocabulary acquisition by example. The builder $A$, when requesting a slab, specifies the linguistic element \w{slab} that he wants to learn. On the other hand, the teacher $B$ associates the right stone to the word \w{slab} and hands it over to $A$. This way, the fibre of $A$ over the word \w{slab} is no longer empty since $A$ has incorporated as part of its meaning the stone received by $B$.
\end{example}
Crucially, in the construction of~\autoref{def:vocabulary acquisition by example} we suppose that $q$ has no meaning at all for the word $L$. What if this is not the case? More generally, we can model a vocabulary acquisition between two speakers that share some prior knowledge of $L$ as a \emph{pushout}, but this will render necessary the introduction of some compatibility conditions. Indeed, consider the following diagram:

\[\vcenter{\xymatrix{
  \clE^q_L \xpo \ar[r]^{u^!}  \ar@{_{(}->}[d] & S \ar[d] \ar@{^{(}->}[r] & \clE^p_L \ar[d]\\
  \clE^q \ar[r] \ar@/_1pc/[dr]_{\framew{q}}& \clF_L \ar@{.>}[d] & \clE^p \ar@/^1pc/[dl]^{\framew{p}}\\
  & \{L\} &
  }}
\]
The arrow marked as $u^!$ always exists whenever $\clE^q_L$ is the empty set because of the universal property of initial objects. Whenever $\clE^q_L$ is not empty, the arrow $u^!$ will have to be explicitly instantiated. This must be understood as: the meaning that $q$ has for $L$ must be compatible with the subset $S$ that constitutes the meaning of the example for $p$. In simpler words, the example $p$ is making \emph{must make sense for $q$}.

\section{Explanations}\label{sec: explanations}
Our task for the remainder of this work is about modelling \emph{vocabulary acquisition by paraphrasis}, which denotes the task of explaining a word by describing it with language, as it happens in a dictionary. To do so, we first need to model what an explanation is.

In stark contrast with works that are exclusively based on syntax, \fiblang can describe linguistic constructions that have meaning for a given speaker \emph{despite their ungrammaticality}. For example, it is a fact of life that often one can understand the meaning of unsound sentences, such as in \emph{`I hungry now'}: this is because every proficient speaker can \emph{interpolate} what is missing in the message they receive by analysing its context, in order to build a grammatical sentence. For that matter, this is precisely how Wittgenstein's language games unravel in disambiguating a sentence.

On the other hand, there are perfectly grammatical sentences, such as the famous \emph{`Dogs dogs dog dog dogs'} (cf. \cite{barton1987}) that are grammatical (since \emph{`dog'} is both a verb and a noun in English) but have no meaning when they are translated to any other speaker out of context. This tension stems from the fact that \emph{acceptability} --i.e. the fact that a sentence has a meaning and \emph{grammaticality} --the fact that a sentence is formed in observance of some generation rules do not fully overlap (cf. \cite{leivada2020}).

As such, leveraging a grammar-based approach to infer meaning in semantics, as in the case of DisCoCat, is going to miss something --an important part, we say. By contrast, the fibrational approach of \fiblang allows more fine-grained bookkeeping: grammaticality is completely encapsulated in the category $\clL$ modelling language, whereas acceptability comes into play in the following definitions:
\begin{definition}[Finite category, finite diagram]
    A \emph{finite category} is a category $\clA$ having a finite set of morphisms. A \emph{finite diagram} valued in $\clL$ is a functor $\clA \to \clL$ whose domain is a finite category.
\end{definition}
\begin{definition}[Explanation]\label{def:explanation}
  Consider a speaker $\var[\framew{p}]{\clE^p}{\clL}$. Fix moreover an object $L$ of $\clL$. An \emph{explanation for $L$ according $p$} is a finite diagram $D_L : \clA \to \clL$ such that the limit $\hat L$ of the diagram
  \[    \xymatrix{
    \clA \ar[r]^{D_L} & \clL \ar[r]^{\GroI{\framew{p}}} & \Set
    }
  \]
  is a subset of the fibre $\clE^p_L$ (Here $\GroI{-}$ is the functor of~\autoref{thm: fibrations functors equivalence}). If $\hat L = \clE^p_L$, we call the explanation \emph{exact}.
\end{definition}
Here, the functor $D_L$ is picking a collection of linguistic elements in the language intending to describe $L$. %For instance, it could select the sentence in~\autoref{eq: pregroup example} of~\autoref{ex: pregroup example}.
The finiteness requirement for $\clA$ stems from the obvious fact that a linguistic sentence is always made of a finite number of words. The linguistic elements picked by $D_L$ are then sent to the sets of sections sitting over them under $\framew{p}$, in accordance with~\autoref{thm: fibrations functors equivalence}. In postulating that $p$ `knows' how to make sense of some given complex concept $L \in \clL$ by breaking it down into some atomic constituents,
it becomes reasonable to assume that the combination of these atomic meanings -- that is, the limit $\hat L$ -- must itself be a concept representing $L$, and thus be a subset of the fibre over it. In light of this interpretation, an exact explanation is a selection of linguistic elements describing the fibre over $L$ completely, that is, an explanation that more than any other, conveys exactly all the nuances that the meaning of $L$ can assume according to the speaker.
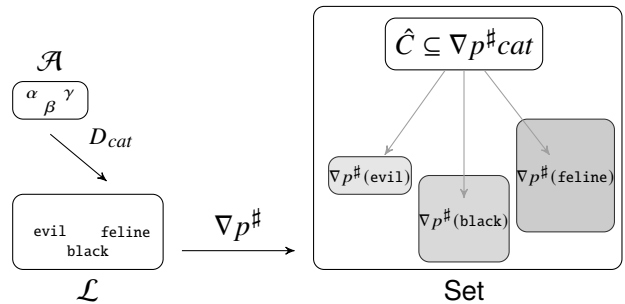
\begin{wrapfigure}{r}{0.65\textwidth}
%\begin{figure}[h]
    \begin{center}
      \begin{tikzpicture}[>=stealth']
        \begin{scope}[xshift=3cm, yshift=2cm]
            \draw[rounded corners] (0,0) rectangle (1,.5) node[pos=.5,right] (inA) {};
            \node[font=\tiny] at (.25,.325) {$\alpha$};
            \node[font=\tiny] at (.5,.15) {$\beta$};
            \node[font=\tiny] at (.75,.325) {$\gamma$};
            \node[above] (Alabel) at (.5,.5) {$\clA$};
        \end{scope}
        \begin{scope}[xshift=3cm]
            \draw[rounded corners] (0,0) rectangle (2,1) node[pos=.5,left] (inL) {};
            \node[below] (Llabel) at (1,0) {$\clL$};

            \node[font=\tiny] at (.5,.5) {$\texttt{evil}$};
            \node[font=\tiny] at (1,.25) {$\texttt{black}$};
            \node[font=\tiny] at (1.5,.5) {$\texttt{feline}$};
        \end{scope}
        \begin{scope}[xshift=7cm]
            \draw[rounded corners] (0,0) rectangle (4.1,3.5) node[pos=.5,right] (Set) {};
            \node[below] (SetLabel) at (2,0) {$\Set$};

            \filldraw[rounded corners, fill=gray!20, draw=black!80] (0.2,1) rectangle (1.3,1.5) node[font=\tiny, pos=.5] (alpha) {$\GroI{\framew{p}} (\texttt{evil})$};
            \filldraw[rounded corners, fill=gray!30, draw=black!80] (1.4,0.1) rectangle (2.6,1.25) node[font=\tiny, pos=.5] (beta) {$\GroI{\framew{p}} (\texttt{black})$};
            \filldraw[rounded corners, fill=gray!40, draw=black!80] (2.7,0.5) rectangle (4,2) node[font=\tiny, pos=.5] (gamma) {$\GroI{\framew{p}} (\texttt{feline})$};

            \node[draw, rectangle,rounded corners] (FL) at  (2,3) {$\hat{C} \subseteq \GroI{\framew{p}} \emph{cat}$};
            %\node[draw, rectangle,rounded corners] (EL) at  (0,3.5) {$\clE_L$};

            \draw[->, thin, gray!75] (FL) -- (alpha);
            \draw[->, thin, gray!75] (FL) -- (beta);
            \draw[->, thin, gray!75] (FL) -- (gamma);
            %\draw[->] (EL) -- (FL) node[font=\tiny, above, pos=.5] {$H_L$};
        \end{scope}

        \draw[->] (3.5,1.8) -- node[font=\footnotesize, above right] {$D_\emph{cat}$} ++(.75,-.625);
        \draw[->] (5.25,0.25) -- node[above] {$\GroI{\framew{p}}$}(6.75,0.25);
      \end{tikzpicture}
    \end{center}
    \caption{An explanation of the word \emph{`cat'} as a \emph{black, evil feline}: The limit of a certain diagram of elements having values in the fibres over concepts like \emph{`black'}, \emph{`evil'}, and \emph{`feline'} is required to be in the fiber over \emph{`cat'}.}
    \label{fig: evilcat}
  \end{wrapfigure}
\begin{remark}[Acceptability vs. grammaticality]
  As we remarked at the beginning of this section, a remarkable feature of~\autoref{def:explanation} is that explanations can be ungrammatical, as they evaluate the acceptability and not grammaticality. For instance, if $\clL$ is a pregroup \cite{093fd03b7dc604c64124630f20a8f01231253397}, the functor $D_L$ can specify a bunch of elements that possibly do not reduce to a sentence type.
\end{remark}

\begin{remark}[On the nature of explanations]\label{rem: on the nature of explanations}
  Explanations are engineered to be far from unique: there may be many functors $D_L$, with different domains, satisfying the property of~\autoref{def:explanation}. This conforms to the idea that the same concept $L$ can be explained in different ways (the same object can be the limit of many different diagrams) by different people, at different moments in time, in different cultures, in different communication settings.

  Interestingly, for every object $L \in \clL$, there is always an exact explanation for it by choosing $\clA$ to be the terminal category and $D_L$ to be the functor picking $L \in \clL$; This explanation is \emph{tautological}, as it affirms in essence that the meaning of the word \emph{`cat'} is \emph{`cat'}.
\end{remark}

\section{Vocabulary acquisition by paraphrasis}\label{sec: vocabulary acquisition by paraphrasis}
As remarked in the opening of~\autoref{sec: explanations}, by \emph{vocabulary acquisition by paraphrasis} we denote the mechanism by which we give meaning to things using linguistic explanation. This is the method of language acquisition commonly used in books and papers like this one, in conversations with a blind man about cathedrals, and pretty much everywhere $p$ needs to convey the meaning of something they know to a $q$ who does not.

\begin{example}\label{ex: vocabulary acquisition by paraprhasis}
  Consider the following dialogue:
  \begin{dialogue}\label{ex:dialogue2}
    \speak{Alice} I adopted a cat!
    \speak{Bob} A what?
    \speak{Alice} You know, a cat: one of those felines.
    \speak{Bob} Oh, you mean, like a tiger?
    \speak{Alice} No: a cat is smaller and it comes in various colours, not only stripes. Mine is black.
    \speak{Bob} Oh, maybe I see. A cat is like\dots a lynx.
    \speak{Alice} Well, almost; black cats are cursed.
    \speak{Bob} Ah, now I see.
  \end{dialogue}
Bob still has probably never seen a real cat when the dialogue ends, and he would have trouble recognizing one. Still, he can get a rough idea of `a cat' by mixing concepts for which he already had a model of meaning.
\end{example}

Mathematically, what we would like to do is proceed as we did in~\autoref{def:vocabulary acquisition by example}, by incorporating an explanation as in~\autoref{def:explanation} into the fibration of the learning speaker. In short, our strategy is the following:
\begin{itemize}
  \item $p$ has an explanation $D_L: \clA \to \clL$ of some word $L$;
  \item $p$ shares $\clA$ with $q$. This represents the act of $p$ uttering the explanation to $q$;
  \item $q$ computes the limit of the diagram $\clA^\op \xrightarrow{D_L^\op} \clL^\op \xrightarrow{\nabla{\framew{q}}} \Set$;
  \item $q$ includes this limit in their own fibre over $L$.
\end{itemize}
Difficulties arise in the last point. The problem we face is that an explanation is a limit, and as such, a particular kind of cone in $\Set$ is composed not only of objects but also of morphisms. Unfortunately, including the whole cone in the total category of a fibration sometimes entails the impossibility of defining the functor $T: \clF^q \to \clL$ as we did in~\autoref{def:vocabulary acquisition by example}.
\begin{example}\label{ex: pregroup no functor}
  To see a practical example, consider the explanation as in~\autoref{fig: evilcat}. Here, the cone legs are morphisms in $\Set$ connecting concepts signifying different \emph{nouns} (for example, from \emph{`cat'} to \emph{`feline'}). If $\clL$ is a Lambek pregroup, the only morphisms are reductions, so there are no morphisms between nouns in $\clL$. Thus in including the whole cone in $\clF^q$, we could not define the functor $T$ on the cone legs. This is related to our previous considerations about grammaticality and acceptability: pregroups only represent grammatical connections and fail to see conceptual relatedness. Consequently, although pregroup grammars have been enriched to identify differently worded sentences in \cite{coecke2021grammar}, the added relations are exquisitely grammatical and cannot account for context.
\end{example}
One possible solution to this problem would be to consider only the limit itself and proceed as in~\autoref{def:vocabulary acquisition by example}. However, in doing so, we would miss the big opportunity of adding meaning while being mindful of the context in which this meaning lives.

Another --more stimulating-- solution considers Wittgenstein's approach to language, which we briefly summarized in the introduction, and leverages the interplay between grammar and semantics. The main idea here is to use semantics data - i.e. the limiting cone - to enrich the grammar with new morphisms. Going back to~\autoref{ex: pregroup no functor}, this means adding new morphisms to the pregroup grammar $\clL$ that do not represent reductions but some sort of `semantic connection' between words.

\def\Quiv{\cate{Quiv}}
\def\Ob#1{{#1}_0}
\def\Id#1{\text{Id}_{#1}}

To add new morphisms to a category, we highlight the following procedure. Recall that there is an adjunction
\[(\firstblank)^\delta :\xymatrix{\Set \ar@<4pt>[r] & \ar@<4pt>[l] \Quiv} : (\firstblank)_0\] sending every quiver $Q$ to the set $Q_0$ of its vertexes, and every set $X$ to the `discrete quiver' $X^\delta$ with no edges, and an adjunction
\[ F :\xymatrix{\Quiv \ar@<4pt>[r] & \ar@<4pt>[l] \Cat} : U\]
sending every quiver $Q$ to the free category $FQ$ generated by it, and every category $\clC$ to its underlying quiver $U\clC$.

% \begin{lemma}
%     \color{red} There is an adjunction $\Set \adjunct{F}{U} \Quiv$ sending every quiver to the set of its vertexes and every set to a `discrete quiver' with no edges.
% \end{lemma}

% \begin{lemma}
%    \color{red}  There is an adjunction $\Quiv \adjunct{F'}{U'} \Cat$ sending every quiver to the free category generated by it, and every category to its underlying quiver.
% s\end{lemma}
The following construction is a recipe to add the set of edges $E$ of a quiver $Q$ to a (small) category $\clC$ and form a category out of it when $\clC$ and $Q$ have the same set of vertices. Regarding a category as a monad in the bicategory of spans, the following \autoref{shit_construction} consists of a particular instance of a free-monoid construction in a decent enough monoidal category (cf. \cite{Lack2010monoids}).

\begin{construction}[FP construction]\label{shit_construction}
Let $\clC$ be a category, and $Q: E \rightrightarrows \Ob{C}$ a quiver over the same set $C_0$ of objects of $\clC$.\footnote{The notation is slightly overloaded here because we denote $Q_0$ the set of vertices of a quiver \emph{and} $C_0$ the set of objects of the underlying quiver of a category $\clC$. This confusion is harmless.}
\begin{itemize}
    \item consider the underlying quiver $U\clC$ of $\clC$;
    \item compute the pushout
    \[\vcenter{\xymatrix{
        \Ob{C}^\delta \xpo \ar[r]\ar[d] & Q \ar[d]\\
        U\clC \ar[r] & U\clC +_{\Ob{C}^\delta} Q
    }}\]
    in the category $\Quiv$ of quivers;%, of the two inclusions $F\Ob{C} \hookrightarrow U'\clC$ and $F\Ob{C} \hookrightarrow Q$, $U\clC +_{F\Ob{C}} Q$;
    \item Applying the free category functor $F : \Quiv \to \Cat$, the square remains a pushout, so we have a pushout
    \[\vcenter{\xymatrix{
        F\Ob{C}^\delta \xpo \ar[r]\ar[d] & FQ \ar[d]\\
        FU\clC \ar[r] & FU\clC +_{F\Ob{C}^\delta} FQ
    }}\]
    in $\Cat$;
    % to it and use the fact that left adjoints commute with colimits: $F'(U\clC +_{F\Ob{C}} Q) \simeq (F'\circ U')\clC \sqcup_{(F' \circ F)\Ob{C}} F'Q$;
    \item compose with the counit $\epsilon$ of the adjunction $F\dashv U$ above:
    \[
        \xymatrix{FU\clC +_{F\Ob{C}^\delta} FQ \ar[r] & \clC +_{F\Ob{C}^\delta} FQ}%(F'\circ U')\clC \sqcup_{(F' \circ F)\Ob{C}} F'Q \xrightarrow{\epsilon' \sqcup_{(F' \circ F)\Ob{C}} \Id{F'Q}} \clC \sqcup_{(F' \circ F)\Ob{C}} F'Q
    \]
\end{itemize}
The category $\clC\wr Q$ is the \emph{FP collage} of $Q,\clC$.
\end{construction}
\begin{figure}[h]
    \begin{center}
        \begin{tikzpicture}
            \draw (0,0) rectangle (2,1.15) node (bx1) {};
              \node (C0) at (-0.3,0.5) {$C_0^\delta$};
              \fill (.5,.5) circle (1pt) node[font = \tiny, left] {$A$};
              \fill (1.5,.75) circle (1pt) node[font = \tiny, right] {$C$};
              \fill (.8,.7) circle (1pt) node[font = \tiny, above] {$L$};
              \fill (1,.25) circle (1pt) node[font = \tiny, right] {$B$};
              \fill (1.75,.25) circle (1pt);
          \begin{scope}[xshift=3cm, yshift=1.5cm]
            \draw (0,0) node (bx2) {} rectangle (2,1.15);
              \node (C) at (-0.4,0.75) {$UC$};
              \fill (.5,.5) circle (1pt) node[font = \tiny, left] {$A$};
              \fill (1.5,.75) circle (1pt) node[font = \tiny, right] {$C$};
              \fill (.8,.7) circle (1pt) node[font = \tiny, above] {$L$};
              \fill (1,.25) circle (1pt) node[font = \tiny, right] {$B$};
              \fill (1.75,.25) circle (1pt);

              \draw[->, >=stealth, shorten <=2pt, shorten >=2pt] (.5,.5) -- (1,.25);
              \draw[->, >=stealth, shorten <=2pt, shorten >=2pt] (1.5,.75) -- (1.75,.25);
              \draw[->, >=stealth, shorten <=2pt, shorten >=2pt] (1,.25) -- (1.5,.75);
          \end{scope}
          \begin{scope}[xshift=3cm, yshift=-1.5cm]
            \draw (0,0) rectangle (2,1.15) node (bx3) {};
              \node (Q) at (-0.25,0.25) {$Q$};
              \fill (.5,.5) circle (1pt) node[font = \tiny, left] {$A$};
              \fill (1.5,.75) circle (1pt) node[font = \tiny, right] {$C$};
              \fill (.8,.7) circle (1pt) node[font = \tiny, above] {$L$};
              \fill (1,.25) circle (1pt) node[font = \tiny, right] {$B$};
              \fill (1.75,.25) circle (1pt);

              \draw[red, ->, >=stealth, shorten <=2pt, shorten >=2pt] (.8,.75) -- (1.5,.75);
              \draw[red, ->, >=stealth, shorten <=2pt, shorten >=2pt] (.8,.75) -- (.5,.5);
              \draw[red, ->, >=stealth, shorten <=2pt, shorten >=2pt] (.8,.75) -- (1,.25);
          \end{scope}
          \begin{scope}[xshift=6cm]
            \draw (0,0) node (bx4) {} rectangle (2,1.15);
            \fill (.5,.5) circle (1pt) node[font = \tiny, left] {$A$};
            \fill (1.5,.75) circle (1pt) node[font = \tiny, right] {$C$};
            \fill (.8,.7) circle (1pt) node[font = \tiny, above] {$L$};
            \fill (1,.25) circle (1pt) node[font = \tiny, right] {$B$};
              \fill (1.75,.25) circle (1pt);
            %==
            \draw[red, ->, >=stealth, shorten <=2pt, shorten >=2pt] (.8,.75) -- (1.5,.75);
            \draw[red, ->, >=stealth, shorten <=2pt, shorten >=2pt] (.8,.75) -- (.5,.5);
            \draw[red, ->, >=stealth, shorten <=2pt, shorten >=2pt] (.8,.75) -- (1,.25);
            \draw[->, >=stealth, shorten <=2pt, shorten >=2pt] (.5,.5) -- (1,.25);
            \draw[->, >=stealth, shorten <=2pt, shorten >=2pt] (1.5,.75) -- (1.75,.25);
            \draw[->, >=stealth, shorten <=2pt, shorten >=2pt] (1,.25) -- (1.5,.75);
          \end{scope}
            \node (push) at (10.5,.5) {$FU\clC +_{FC_0^\delta} FQ$};
            \node[below=1cm of push] (wreath) {$\clC\wr Q$};
            \draw[->, >=stealth] (push) -- node[right] {$\epsilon_C +_{FC_0^\delta} 1_{FQ}$} (wreath);
            \draw[->, >=stealth] (bx1) -- (bx2);
            \draw[->, >=stealth, shorten >=1mm, shorten <=1mm] (2,0) -- (3,-.5);
            \draw[->, >=stealth, shorten >=1mm, shorten <=1mm] (5,1.5) -- (6,1);
            \draw[->, >=stealth] (bx3) -- (bx4);
            \draw[->, line join=round, decorate,
                  decoration={
                                zigzag,
                                segment length=4,
                                amplitude=.9,post=lineto,
                                post length=2pt
                              },
                  >=stealth] (8.1,.5) -- node[above] {$F$} (push.west);
            \end{tikzpicture}
    \end{center}
    \caption{A graphical description of $\clC\wr Q$, with the pushout of quivers made explicit.}
\end{figure}
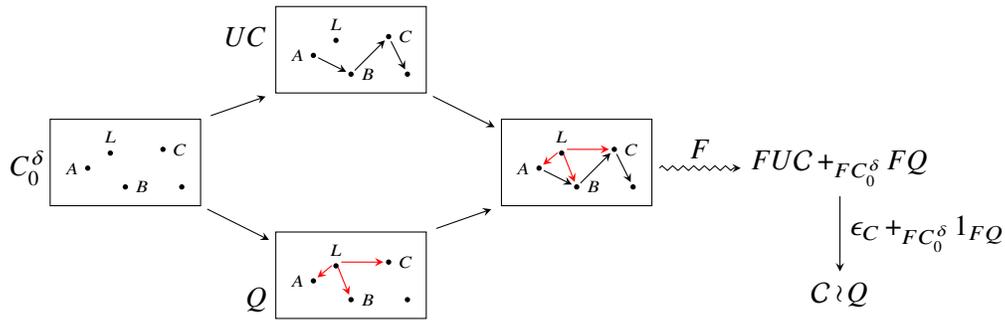
\begin{remark}
Note that there exists a canonical functor
\[\xymatrix{
    K : \clC \ar[r] & \clC + FQ \ar[r] & \clC \wr Q
}\] given by the coproduct embedding followed by the projection on the quotient realising the pushout $\clC\wr Q$; note that by construction this functor is the identity on objects, so it is induced in a canonical way by a monad $\fkq : \xymatrix{\clC \ar|-*=0@{|}[r] & \clC}$ on $\clC$ in the category of profunctors, and $K$ corresponds to the free functor into the Kleisli category of $\fkq$.
\end{remark}
\begin{remark}
A more general construction for $\clC\wr Q$ is then the following: fix a monad $\fkq$ as above, and consider its Kleisli object $\bar\clC$; the free part of the Kleisli adjunction yields an identity on objects functor $\clC \to \bar\clC$. Given $\fkq$, it can be highly complicated to describe its Kleisli category; what makes this construction combinatorially tamer is that the structure we are adding through the quiver $Q$ is free.
\end{remark}

Now we finally have all the needed tools to define vocabulary acquisition via paraphrasis satisfactorily.

\begin{definition}[Vocabulary acquisition by paraphrasis]\label{def: vocabulary acquisition by paraphrasis}
    Consider two speakers $\var[\framew{p}]{\clE^p}{\clL}$ and $\var[\framew{q}]{\clE^q}{\clL}$, which we will call \emph{teacher} and \emph{learner}, respectively.

    Suppose that, for some $L \in \clL$ --called \emph{the linguistic element to learn}-- we have that $\clE^p_L \neq \emptyset$ and $\clE^q_L = \emptyset$. Let $D_L : \clA \to \clL$ be an explanation of $L$ according to $p$. Define the category $\clF^q$ as $\clL \wr Q$, where $Q$ is the quiver obtained as follows:
    \begin{itemize}
      \item the vertices are the same of $\clL$;
      \item there is an edge $L\to L'$ for each limiting cone leg $\lim \big(\clA^\op \to \clL^\op \xto{\GroI \framew{q}} \Set\big) \to \nabla\framew{q} L'$.
    \end{itemize}

    Now define a functor $T: \clF^q \to \Set$ by mapping $L$ to $\hat{L}^q$ and every other $L'$ to $\nabla\framew{q}(L')$. On morphisms, $T$ agrees with $\nabla\framew{q}$ wherever the latter is defined and maps the newly added edges of $Q$ to the legs of the limit $\hat{L}^q$.
  \end{definition}
  Using~\autoref{thm: fibrations functors equivalence}, the new fibration modelling $q$ after learning $L$ is $\int T$.

  \medskip
This definition is a bit terse and needs some unpacking, so let us piggyback to~\autoref{ex: vocabulary acquisition by paraprhasis}. In our situation, $p$ knows at least partly what a \emph{`cat'} (the object $\catIcon$) is, because the fibre $\clE^p_{\catIconSmall}$ is not empty. We postulate that $p$ can explain this concept in words, that is, $p$ can define a functor $D_{\catIconSmall}$ picking a bunch of linguistic elements in the language that mean \emph{`cat'} to $p$ since, by~\autoref{def:explanation}, the limit $\hat{\catIcon}^p$ of $\GroI{\framew{p}} \circ D_{\catIconSmall}$ is a subset of $\GroI{\framew{p}} \catIcon$, which corresponds exactly to $\clE^p_{\catIconSmall}$ via~\autoref{thm: fibrations functors equivalence}.

The teacher $p$ would like to `transmit' $\hat{\catIcon}^p$ to $q$, but this is not possible since, unless we postulate either of them is from the planet Vulcan, the only way $p$ and $q$ have to communicate is through the language $\clL$. Still, $p$ can utter the explanation and thus share the functor $D_{\catIconSmall}$ with $q$. Notice how, at this stage, $D_{\catIconSmall}$ is most likely \emph{not} an explanation for $q$ as $\GroI{\framew{q}} \catIcon = \clE^q_{\catIconSmall}$ is empty by definition.

In any case, $q$ can calculate the limit $\hat{\catIcon}^q$ of $\GroI{\framew{q}} \circ D_{\catIconSmall}$. The limits $\hat{\catIcon}^q$ and $\hat{\catIcon}^p$ will, in general, be different, as the same explanation makes sense to different speakers in different ways. Notice that whereas $\hat{\catIcon}^p$ is always non empty as $D_{\catIconSmall}$ is defined to be an explanation for $p$, $\hat{\catIcon}^q$ can be empty. This happens when $q$ is not able to successfully combine the meanings of the words in the explanation $D_{\catIconSmall}$: the explanation does not make sense to $q$. Interestingly, this is the case for the tautological explanation of~\autoref{rem: on the nature of explanations}, which captures perfectly the meaning of $\catIcon$ for $p$, but means absolutely nothing to $q$.

Whenever $\hat{\catIcon}^q$ is not empty, this is nothing more than a combination of concepts that $q$ already knows, as illustrated in~\autoref{fig: evilcat}. $q$ includes this composition of concepts in the fibre over $L$, while the morphisms from the limit to its atomic constituents are included as morphisms between the fibres. In this procedure, the underlying language category for $q$ changes, as we now have a fibration over the category $\clF^q$, which is obtained by adding new morphisms to the language category $\clL$. This is not a bug but a feature: in learning the meaning of a new concept, the speaker $q$ also learns new ways to turn words and sentences into others.

\begin{remark}\label{rem: slab}
  The language game in~\autoref{exa:language game} is taken up and debated later in \cite{wittgenstein2010philosophical} with a question: should we interpret the request `\emph{Slab!}' as a word or a sentence? In the latter case, the sentence `\emph{Slab!}' should be understood as a shortening of the sentence `\emph{Bring me a slab!}'.
  Concerning this sentence, we can consider the roles as inverted with respect to~\autoref{exa:language game}: the builder makes a request, and the assistant can use~\autoref{def: vocabulary acquisition by paraphrasis} to build a meaning for the sentence `\emph{Bring me a slab!}'. Hopefully, this meaning will include things contextually relevant to the situation, allowing the assistant to identify the slab and fulfil the builder's request correctly.

  Yet, this does not answer Wittgenstein's question, namely how we manage to go from `\emph{Slab!}' to `\emph{Bring me a slab!}'. Wittgenstein argues that the reason behind this semantic identification is the fact that these two sentences admit the same \emph{contextual use}: there is no need for an explicit explanation, as context is directly responsible of disambiguating this sentence.
\end{remark}

  ...But what is context mathematically? This is quite a thorny question. One possible solution to model context in our framework is using stronger models for the language category $\clL$. In many examples, $\clL$ is taken to be a pregroup, which corresponds to a context-free grammar in the sense of Chomsky \cite{Buszkowski2001}. We could instead make use of models of grammar that are \emph{context-sensitive}, allowing for a finer degree of context management. In theory, we could have grammars where `\emph{Slab!}' can be reduced to `\emph{Bring me a slab!}' in a given context, and proceed as specified in~\autoref{rem: slab}.

  Yet, if there is something that~\autoref{def: vocabulary acquisition by paraphrasis} taught us, it is that the process of acquiring vocabulary can result in enriching language with semantic meaning. Following the same idea, we could start with a model of $\clL$ that is context-free, such as a pregroup, and gradually adding `context-sensitive' morphisms that we borrow from the semantics, exactly as in~\autoref{def: vocabulary acquisition by paraphrasis}. Going back to~\autoref{rem: slab}, the reduction from `\emph{Slab!}' to `\emph{Bring me a slab!}' would be added to $\clL$ as a result of a language game previously played.

  This last consideration hints at the fact that we shall be able to define disambiguation, and more in general communication, fibrationally. This is a broad topic, and a current matter of investigation.
\section{Conclusion and future work}
In this work, we used the framework provided by FibLang to take a first stab at describing vocabulary acquisition mathematically. In particular, we defined the concept of \emph{explanation} to define, in turn, vocabulary acquisition by paraphrasis. A clear direction of future work is to use~\autoref{def:explanation} as the building block for a more general theory of communication: intuitively, speakers communicate in a game-theoretic fashion, exchanging explanations and using them to build meanings until they reach some kind of fixed point, at which providing further explanations does not result in building further meaning. This broader picture would fully capture Wittgenstein's ideas regarding language games, and we consider it an ambitious goal.
With respect to this, a conjecture we are currently working on is to show how the bidirectional, game-theoretic nature of interactions between two speakers naturally suggests the use of self-dual categorical structures.

On a more practical standpoint, we would like to investigate a possible replacement of the limit in~\autoref{def:explanation} and in~\autoref{def: vocabulary acquisition by paraphrasis} with something more linguistically sound: it is true that limits are universal constructions in category theory, and, as such, are mathematically well-behaved. Yet, in strictly applied contexts it may be useful to experiment with alternative definitions: for instance, considering a model of meaning where concepts in the fibres are images may be a sensible choice to experiment with machine-learning algorithms that merge images together \cite{Gatys2016}.

\appendix

\printbibliography

\end{document}